\documentclass[11pt]{amsart}

\newtheorem{theorem}{Theorem} %[section]
\newtheorem{prop}[theorem]{Proposition}
\newtheorem{lemma}[theorem]{Lemma}
\newtheorem{cor}[theorem]{Corollary}
 
 \usepackage{amsfonts}
\usepackage{amsopn}
\usepackage{epsfig}
\usepackage{color}
\usepackage{amsmath}
\usepackage{graphicx}
\usepackage{amssymb}

\newcommand\C{{\bf C}}
\newcommand\Chat { {\hat{\C}} } 

\renewcommand\P{{\bf P}}
\newcommand\R{{\bf R}}
\newcommand\Z{{\bf Z}}
\newcommand\D{{\bf D}}
\newcommand\del{\partial} 
\newcommand\delbar{\bar{\del}}

\newcommand\eps{\varepsilon}

\renewcommand\phi{\varphi}

\newcommand\<{\langle} %innerproduct brackets 
\renewcommand\>{\rangle} %%

\newcommand\iso{\simeq} %isomorphism -~
\renewcommand\O{\mathcal{O}}

\newcommand\scriptC{\mathcal{C}}

\newcommand\PSL{\operatorname{PSL}}

\newcommand\PGL{\operatorname{PGL}}
\renewcommand\H{\operatorname{H}}

\newcommand\Id  {\operatorname{Id}}     %Identity map
\newcommand\supp{\operatorname{supp}}   %support
\newcommand\dist{\operatorname{dist}}  %distance function
\newcommand\tr  {\operatorname{tr}}
\renewcommand\gcd {\operatorname{gcd}} %gcd
\newcommand\Res {\operatorname{Res}} %resultant
\newcommand\Rat  {\operatorname{Rat}}
\newcommand\Poly {\operatorname{Poly}}
\newcommand\Ratbar {\overline{\Rat}}  
  
\newcommand\Polybar{\overline{\Poly}}

\begin{document}

\title[Iteration at the boundary]{Iteration at the 
boundary of the space of rational maps}

\author{Laura DeMarco}
\date{December 20, 2004}

\begin{abstract}
Let $\Rat_d$ denote the space of holomorphic self-maps
of $\P^1$ of degree $d\geq 2$, and $\mu_f$ the measure 
of maximal entropy for $f\in\Rat_d$.  The map of measures
$f\mapsto\mu_f$ is known to be continuous on $\Rat_d$, and
it is shown here to extend continuously to the boundary of
$\Rat_d$ in $\Ratbar_d \iso 
%$\P\H^0(\P^1\times\P^1, \O(d,1))\iso
\P^{2d+1}$,
except along a locus $I(d)$ of codimension $d+1$.  The set 
$I(d)$ is also the indeterminacy locus of the iterate 
map $f\mapsto f^n$ for every $n\geq 2$. The limiting measures are given
explicitly, away from $I(d)$.  The
degenerations of rational maps are also 
described in terms of metrics of non-negative curvature on the 
Riemann sphere:  the limits are polyhedral. 
\end{abstract}

\maketitle

\thispagestyle{empty}

\bigskip
For each integer $d\geq 1$, let $\Rat_d$ denote the space of holomorphic maps
$f:\P^1\to\P^1$ of degree $d$ with the topology of uniform convergence.  
Fixing a coordinate system on the 
projective line, each such map can be expressed as a ratio of 
homogeneous 
polynomials $f(z:w) = (P(z,w):Q(z,w))$, where $P$ and $Q$ have no common 
factors and are both of degree $d$.  Parameterizing the space
$\Rat_d$ by the coefficients of $P$ and $Q$, we have 
     $$\Rat_d \iso \P^{2d+1} - V(\Res),$$
where $V(\Res)$ is the hypersurface of polynomial pairs $(P,Q)$ for which
the resultant vanishes.  In particular, $\Rat_d$ is smooth and affine.

In this paper, we aim to describe the possible limiting behavior of a 
sequence of rational maps which diverges in $\Rat_d$, for each 
$d\geq 2$, in terms of the 
measures of maximal entropy and corresponding conformal metrics 
on the Riemann sphere.  This is the
first step in describing a natural compactification of this space, or 
a boundary of the moduli space $\Rat_d/\PSL_2\C$, which
is well-behaved under iteration.   A compactification of the moduli space 
has been studied by Milnor \cite{Milnor:quad} and 
Epstein \cite{Epstein:bounded} in degree 2 and Silverman
\cite{Silverman} in all degrees, but 
iteration does not extend continuously
to this boundary, as first seen in \cite{Epstein:bounded}. 
See \cite{D:moduli2} for more details.  

We can associate to every point in 
$\Ratbar_d\iso\P^{2d+1}$ a 
self-map of the Riemann sphere of degree $\leq d$, together with
a finite set of marked points.   
Namely, each $f\in\Ratbar_d$ determines the coefficients for a pair of 
homogeneous polynomials, defining a map on $\P^1$  
away from finitely many {\bf holes}, the shared roots of the pair
of polynomials.  Each hole comes  
with a multiplicity, the {\bf depth} of the hole.  
We also define a probability measure $\mu_f$ 
for each $f\in\Ratbar_d$.  For $f\in\Rat_d$, we let $\mu_f$ be
the unique measure of maximal entropy for $f$ so that $\mu_f = 
\mu_{f^n}$ for all iterates of $f$.  See 
\cite{Lyubich:entropy},\cite{FLM}, and \cite{Mane:unique}.  For each
$f$ in the complement of $\Rat_d$, the measure $\mu_f$ will be
a countable sum of delta masses, with atoms at the holes of $f$ and 
along their backward orbits.  See Section \ref{definitions} for the definitions.

The {\bf indeterminacy  locus} $I(d)\subset\Ratbar_d$
in the boundary of $\Rat_d$ consists of degree 0 maps such 
that the constant value is also one of the holes.  
The codimension of $I(d)$ is $d+1$.   For each 
$f \not\in I(d)$, we will see in Section \ref{iteration} that $\mu_f
=\mu_{f^n}$ for all iterates of $f$.  
We prove, 

\begin{theorem}  Fix $d\geq 2$, and 
suppose $\{f_k\}$ is a sequence in $\Rat_d$ converging in $\Ratbar_d$
to $f\not\in\Rat_d$. 
\begin{itemize}
\item[(a)] For $f\not\in I(d)$, the measures of maximal
           entropy $\mu_{f_k}$ for $f_k$ converge weakly to $\mu_f$.  
\item[(b)] For $f\in I(d)$, 
           any subsequential limit $\nu$ of the measures 
           $\mu_{f_k}$ must 
           satisfy $\nu(\{c\}) \geq d_c/(d_c+d)$, where 
$c\in\P^1$ is both the constant value of $f$ and a hole 
of depth $d_c\geq 1$.  
\end{itemize}
\end{theorem}

\noindent 
If $f\not\in I(d)$ has a hole at $h\in\P^1$ of depth $d_h\geq 1$, then
$\mu_f(\{h\}) \geq d_h/d$.
In Section \ref{Thm2}, Example 1, we provide  examples in $\Rat_d$ 
for every $d\geq 2$ which realize the 
lower bound of part (b) when $d_c=1$, so that $\nu(\{c\})=1/(d+1)$.  

\medskip\noindent  {\bf The iterate map}.
Theorem 1 is, in part, an extension of a result of Ma\~n\'e which 
states that the measures of maximal entropy vary continuously (in 
the weak topology) over $\Rat_d$ \cite[Thm B]{Mane}.  
The proof of Theorem 1 relies on the study of the 
iterate map  $\Phi_n: \Rat_d \to \Rat_{d^n}$
which sends a rational map $f$ to its $n$-th iterate $f^n$.  The 
iterate map $\Phi_n$ extends to a rational map from $\Ratbar_d$ to 
$\Ratbar_{d^n}$.  We obtain,

\begin{theorem}  For each $d\geq 2$, the following are equivalent:
\begin{itemize}
\item[(i)]  $g\in\Ratbar_d$ is in the indeterminacy locus $I(d)$,
\item[(ii)]  the iterate map $\Phi_n$ is undefined at $g$ for some $n\geq 2$,
\item[(iii)]  the iterate map $\Phi_n$ is undefined at $g$ for
               all $n\geq 2$, and
\item[(iv)]  the map $f\mapsto\mu_f$ is discontinuous at $g$.
\end{itemize}
\end{theorem}

In particular, the map of measures $f\mapsto\mu_f$ extends 
continuously from $\Rat_d$ to a point $g$ in the boundary if and only 
if the iterate map $\Phi_n$ extends continuously to $g$ for some
$n\geq 2$.
The understanding of the iterate map also motivated Theorem 1, 
and from it we obtain the following corollary.

\begin{cor} \label{proper}
The iterate map $\Phi_n: \Rat_d \to \Rat_{d^n}$ given by
$f\mapsto f^n$ is proper for all $n$ and $d\geq 2$.  
\end{cor}

\proof
The measure of maximal entropy $\mu_f$ for $f\in\Rat_d$ is always
non-atomic, and the map $f\mapsto\mu_f$ is continuous (with the 
topology of weak convergence on the space of probability measures on the
Riemann sphere) on $\Rat_d$ \cite[Thm B]{Mane}.  Now,
suppose $\{f_k\}$ diverges in $\Rat_d$.  There exists a subsequence of the 
maximal measures $\mu_{f_k}$ 
which converges weakly to a measure $\nu$.  By 
Theorem 1, $\nu$ has atoms.  Recall that for $f\in\Rat_d$,
the measure $\mu_f$ is also the 
measure of maximal entropy for all the iterates of $f$.   
If for some $n$ the sequence of iterates $\{f_k^n\}$ converges 
in $\Rat_{d^n}$,
the measures $\mu_{f_k}$ would have to converge to a non-atomic 
measure. \qed

\medskip
Note that properness does {\em not} hold in degree 1, 
where $\Rat_1 \iso \PGL_2\C$, since there are 
unbounded families of elliptic M\"obius maps of finite order.
For example, the family $f_k(z) = k/z$ diverges in $\Rat_1$
as $k\to \infty$, but $f_k^2(z) = z$ for all $k$.  
Properness of the iterate map for degrees $>1$ should be 
intuitively obvious, 
but analyzing its behavior near $I(d)$ is somewhat delicate.

%Currently under investigation is the resolution of the indeterminacy of
%the iterate map for each $n$.  It would be interesting to know when,
%for example, a sequence
%limiting in $I(d)$ has an iterate converging to a map of positive
%degree.  Even in the space of polynomials, interesting examples appear
%which can give information about how the multipliers at periodic points
%behave.  

\medskip\noindent {\bf Dynamics at the boundary}.  
In Section \ref{dynamics}, we describe the dynamics of a  map 
in the boundary of $\Rat_d$.  In particular, for
$f\not\in I(d)$, we can define the escape rate function $G_F$
of a homogenization $F:\C^2\to\C^2$ of $f$, 
and we show that it satisfies $dd^c G_F = \pi^*\mu_f$,
just as for non-degenerate rational maps (see 
\cite[Thm 4.1]{Hubbard:Papadopol}).   
In analogy with results for rational maps,  the
measure $\mu_f$ is also the limit of pull-backs by the iterates
of $f$ of any probability
measure on $\P^1$.

\medskip\noindent {\bf Metric convergence}. 
Every rational map determines a conformal 
metric on $\P^1$ (unique up to scale) with 
non-negative distributional curvature equal to the measure of maximal 
entropy (see Section \ref{metrics}).  The sphere with this metric can 
be realized as the intrinsic metric of a convex surface in $\R^3$ 
by a theorem of Alexandrov \cite[VII \S 7]{Alexandrov}.   
Each $f\in \Ratbar_d - \Rat_d$ determines a conformal metric with 
singularities, defined by the measure 
$\mu_f$.  Theorem 1 together with 
\cite[Thm 7.3.1]{Reshetnyak} implies,

\begin{cor} \label{metricprop}
Suppose $f_k$ in $\Rat_d$ converges in $\Ratbar_d$ to $f\not\in I(d)$. 
Then the spheres with associated metrics have a convex polyhedral 
limit with distributional curvature $4\pi\mu_f$.   
\end{cor}

\noindent
The (countably many) cone points in the limiting metric of 
$f\not\in I(d)$ have cone angles given by 
$2\pi - 4\pi\mu_f(\{p\})$ for each $p\in\P^1$, 
where a non-positive cone angle means that the point is
at infinite distance from all others.    
The metric convergence is uniform away 
from the infinite ends.  

The subspace of polynomials $\Poly_d\subset\Rat_d$ is very interesting 
by itself as the limiting measures (away from $I(d)$) are supported
at $\leq d$ points.  
%should be always atomic with {\em finitely many} points in the support.  
See Section \ref{examples}. 
The metrics associated to such measures
are polyhedral with finitely many vertices.  Thus, a boundary of the
moduli space of polynomials can be described in part by the geometry
of convex polyhedra (as in e.g. \cite{Thurston:shapes}).  
%This perspective will be pursued further in a sequel.

\medskip\noindent{\bf Outline.}  In Section 1, we fix 
notation and define the probability measure $\mu_f$ for each 
$f\in\Ratbar_d$.  Section 2 is devoted to a study of the indeterminacy 
of the iterate map at the boundary of $\Rat_d$.  In Section 3, we 
study the dynamics of a map in the boundary, and we show the existence
of the escape rate function.  Theorem 1 is proved in Section 4, and 
Theorem 2 is proved in Section 5.  The Alexandrov geometry of
rational maps and Corollary \ref{metricprop} are described in Section
6.  We conclude with some further examples in Section 7.  

I would like to thank Adam Epstein for all of his helpful and creative 
comments and whose paper \cite{Epstein:bounded} prompted my study of 
the iterate map.   I'd also like to thank Ilia Binder, Mario Bonk, and 
Curt McMullen for their suggestions.

%%%%%%%%%

\section{Definitions and notation} \label{definitions}

There is a natural map 
  $$\Rat_d \hookrightarrow 
\Gamma(d,1) = \P\H^0(\P^1\times\P^1, \O(d,1)),$$
of the space of rational maps into the projectivized space 
of global sections of the line bundle 
$\O(d,1)$ on $\P^1\times\P^1$.  The graph of
$f:\P^1\to\P^1$ defines the zero locus of a section. 
In coordinates, the rational map defined by 
  $$f(z:w) = (P(z,w):Q(z,w))$$
is sent to the section represented by 
  $$\P^1\times\P^1\ni ((z:w),(x:y)) \mapsto
        yP(z,w) - xQ(z,w).$$
Consequently, we have 
  $$\Ratbar_d = \Gamma(d,1) \iso \P^{2d+1}.$$

Every $f\in\Ratbar_d$ determines the coefficients for a pair 
of homogeneous polynomials, and we write
  $$f = (P:Q) = (H p: H q) = H\phi_f,$$
where $H = \gcd(P,Q)$ is a homogeneous polynomial and $\phi_f
= (p:q)$ is a rational map of degree $\leq d$.  A zero of 
$H$ in  $\P^1$ is said to be a {\bf hole} of $f$, and
the multiplicity of such a zero is its {\bf depth}.  The holes
can be interpreted as punctures in the domain of definition
of $f$ as a map from $\P^1$ to itself (though as singularities,
they are removable).  When $f\in\Ratbar_d$ has holes, 
it is said to be {\bf degenerate}.  
The graph of a degenerate $f = H\phi_f \in\Ratbar_d$ is given by
 $$\Gamma_f = \{(p,\phi_f(p)) \in\P^1\times\P^1\} \cup 
                 \{ h\times\P^1: H(h) = 0\},$$
and has vertical components (counted with multiplicity) over the holes
of $f$.   We see that $\Gamma_f$ defines a holomorphic correspondence
in $\P^1\times\P^1$, and a hole can also be interpreted as a point
which is mapped by $f$ over the whole of $\P^1$.  

The {\bf indeterminacy locus}
$I(d)\subset\Ratbar_d$ is the set of degenerate maps $f=H\phi_f$ 
for which $\phi_f$ is constant and this constant value is one
of the holes of $f$;  that is, $f$ has the form 
$f(z:w) = (aH(z,w):bH(z,w))$ for some $(a:b)\in\P^1$ 
with $P(a,b)=0$, and therefore, $I(d)$ is given by
$$I(d) = \{ f = H \phi_f :  
\deg\phi_f = 0 \mbox{ and } \phi_f^*H \equiv 0\}.$$
A simple dimension count shows that 
$I(d)$ has codimension $d+1$.  In fact, 
the locus $I(d)$ is isomorphic to $\P^1\times\P^{d-1}$
by sending $f = H\phi_f \in I(d)$ with $\phi_f\equiv (a:b)$  
to the pair $((a:b), H(z,w)/(bz-aw))$.  See Figure \ref{figure1}.

\bigskip
\begin{figure}[htbp]
\begin{center}
\input{Ratbar.pstex_t}
\end{center}
\caption{Graphs in $\P^1\times\P^1$ of $f\in \Ratbar_3$:  (a) $f\in\Rat_3$, 
(b) $f = H\phi_f\in\del\Rat_3$ with $\deg\phi_f = 2$, (c) $f\in I(3)$.}
\label{figure1}
\end{figure}
%\bigskip

\medskip\noindent{\bf Example.}  For $d=1$, we have
$\Ratbar_1\iso\P^3$, the space of all non-zero two-by-two
matrices $M$ up to scale.  The indeterminacy locus is  
  $$I(1) = \{M: \tr M = \det M = 0\}.$$
Indeed, if $M\in I(1)$, by a change of coordinates, 
we can assume that $M$ is the 
constant map infinity on $\P^1$ with one hole at infinity.
In coordinates, $M(z:w)=(w:0)$, or rather, $M$ is a matrix
with one non-zero entry off the diagonal and zeroes elsewhere.
Up to conjugacy, these are precisely the matrices with 
vanishing trace and determinant.

\medskip\noindent
{\bf Probability measures}.
We define a probability measure 
$\mu_f$ on $\P^1$ for each point $f\in\Ratbar_d$.  
For $f\in \Rat_d$, a point $a\in\P^1$
is {\em exceptional} if its grand orbit $\bigcup_{n\in\Z} \{f^n(a)\}$
is finite.  By Montel's Theorem, a rational
map can have at most two exceptional points.  
Let $\mu_f$ be 
the unique measure of maximal entropy for $f$, given by the weak limit,
  $$\mu_f = \lim_{n\to\infty}  \frac{1}{d^n} 
               \sum_{f^n(z)=a}   \delta_z,$$
for any non-exceptional point $a\in\P^1$  
\cite{Lyubich:entropy},\cite{FLM},\cite{Mane:unique}.  The measure
$\mu_f$ has no atoms, and its support equals the Julia set of $f$; 
it is also the unique measure of maximal entropy for every 
iterate $f^n$.  

If $f = H\phi_f$ is degenerate and $\phi_f$ 
is non-constant, we define an atomic measure, 
 $$\mu_f := \sum_{n=0}^\infty \frac{1}{d^{n+1}} 
                \sum_i \sum_{\phi_f^n(z)=h_i} \delta_z,$$
where the holes $h_i$ and all preimages by $\phi_f$ are counted
with multiplicity.  Note that if the hole $h_i$ has depth $d_i$, then 
$\mu_f(\{h_i\})\geq d_i/d$.  Furthermore, since $\deg \phi_f = 
d - \sum_i d_i$, the total measure $\mu_f(\P^1)$ is 1.
If $\phi_f$ is constant, then the depths 
of the holes sum to $d$, and we set 
 $$\mu_f = \frac{1}{d} \sum_i  \delta_{h_i},$$
where again the holes $h_i$ are counted with multiplicity.  

We will see in the following section that for 
every degenerate $f\not\in I(d)$, 
we will have $\mu_f = \mu_{f^n}$ for all iterates of $f$.  

\medskip\noindent
{\bf Example.}  Suppose that $f(z:w) = (P(z,w):w^d)$ 
where $P\not\equiv 0$ is a homogeneous polynomial such that
$P(1,0)=0$.  Then $f$ is degenerate with a hole at $\infty = 
(1:0)$.  The map $\phi_f$ can be identified with a polynomial in $\C$
of degree $< d$,
by choosing local coordinates $z/w$ for $(z:w)\in\P^1$.
Since the backwards orbit of $\infty$ under any polynomial
consists only of $\infty$ itself, we must have 
$\mu_f = \delta_\infty$.

\medskip
Generally, when computing the $\mu_f$-mass of a point for
degenerate $f\in\Ratbar_d$, one needs to count the number 
of times the forward
iterates of the point land in a hole of $f$.   The following lemma 
follows directly from the definition of the measure $\mu_f$.  

\begin{lemma} \label{pointmass}
Let $f= H\phi_f\in\Ratbar_d$ be degenerate with $\deg\phi_f>0$.
For each $a\in\P^1$, we have,
 $$\mu_f(\{a\}) = \frac{1}{d} \sum_{n=0}^\infty 
          \frac{ m(\phi_f^n(a)) d(\phi^n(a))}{d^n},$$
where $d(\phi_f^n(a))$ is the depth of $\phi_f^n(a)$ as a hole of 
$f$ and $m(\phi_f^n(a))$ is the multiplicity of $a$ as a solution
of $\phi_f^n(z) = \phi_f^n(a)$.  
\end{lemma}

The space $M^1(\P^1)$ of probability measures on $\P^1$ is given 
the weak topology.  For what follows, it will be useful to recall that
$M^1(\P^1)$ is metrizable because it is 
a compact subset of the dual space to the separable $C(\P^1)$, the 
continuous functions on $\P^1$.

%%%%%%%%%

\section{Iterating a degenerate map} \label{iteration}

The iterate
map $\Phi_n: \Rat_d\to\Rat_{d^n}$, which sends $f$ to $f^n$, 
is a regular morphism between smooth affine varieties.  It 
extends to a rational map $\Ratbar_d\dashrightarrow\Ratbar_{d^n}$
for all $n\geq 1$.
In this section, we will give a formula for the iterates of 
a degenerate map, where defined, and specify the indeterminacy 
locus of the iterate map $\Phi_n$.  

Let $(a_d,\ldots,a_0,b_d,\ldots,b_0)$ denote the 
homogeneous coordinates on $\Ratbar_d\iso\P^{2d+1}$, where
a point $f = (P:Q)$ has coordinates given by the coefficients
of $P$ and $Q$.     
The $2d^n+2$ coordinate functions which define the iterate 
map $\Phi_n:\Ratbar_d\dashrightarrow\Ratbar_{d^n}$ 
generate a homogeneous ideal $I_n$
in the ring $A={\bf Z}[a_d,\ldots,b_0]$.  
The ideal $I_1 = (a_d,\ldots,b_0)$,  generated by all 
homogeneous monomials of degree 1 in A, is the ideal generated 
by the identity map $\Phi_1$.  

\begin{lemma}  \label{descending}  
In the ring $A$, the ideals $I_n$ are generated by 
homogeneous polynomials of degree $(d^n-1)/(d-1)$ and satisfy 
  $$I_n \subset I_1\cdot I_{n-1}^d $$
for all $n\geq 2$.  In particular, they form a 
descending chain.  
\end{lemma} 

\proof
The affine space $\C^{2d+2}$ parameterizes, by the coefficients,
all pairs $F=(P,Q)$ of degree $d$ homogeneous polynomials 
in two variables.  Such a pair defines a map $F:\C^2\to\C^2$,
and the composition map
 $$\scriptC_{d,e}: \C^{2d+2}\times\C^{2e+2} \to \C^{2de+2}$$
sending $(F,G)$ to the coefficients of $F\circ G$ is bihomogeneous of degree 
$(1,d)$ in the coefficients of $F$ and $G$.  In particular, 
the second iterate $\Phi_2$ is (the projectivization of) the 
restriction of 
$\scriptC_{d,d}$ to the diagonal of $\C^{2d+2}\times\C^{2d+2}$,
and so its coordinate functions are homogeneous of degree $1+d$.  
Thus, the ideal 
$I_2$ in $A$ generated by these coordinate functions of $\Phi_2$
is contained in $I_1^{d+1} = I_1\cdot I_1^d$.  

For the general iterate, of course $F^n = F\circ F^{n-1}$, 
so $\Phi_n$ can be expressed as 
 $$\Phi_n = \scriptC_{d,d^{n-1}} \circ (\Id,\Phi_{n-1}):
               \C^{2d+2} \to \C^{2d^n+2}.$$
Consequently, $\Phi_n$ is homogeneous of degree $1 +
d(\deg\Phi_{n-1})$.  By induction, we have $\deg\Phi_n 
= 1+d+\cdots+d^{n-1} = (d^n-1)/(d-1)$.  The above expression
for $\Phi_n$ and the bihomogeneity of the composition map 
implies that the coordinate functions of $\Phi_n$ must 
lie in the ideal $I_1\cdot I_{n-1}^d$.  
\qed

\medskip
Recall from Section \ref{definitions} that $I(d)\subset\Ratbar_d$
is defined as the locus of degenerate constant maps such that
the constant value is equal to one of the holes. 

\begin{lemma} \label{iterate}
The indeterminacy locus for the iterate map 
$\Phi_n: \Ratbar_d\dashrightarrow\Ratbar_{d^n}$
is $I(d)$ for all $n\geq 2$ and all $d\geq 1$.  
If $f=H \phi \not\in I(d)$ is degenerate, then 
  $$f^n = \left( \prod_{k=0}^{n-1} (\phi^{k*}H)^{d^{n-k-1}} \right) \phi^n.$$
\end{lemma}

\proof
Suppose that $f = (P:Q) = (H p:H q) = H \phi$ is 
degenerate.  The second iterate of $f$ has the form,
\begin{eqnarray*}
  f\circ f  &=&  (P(P,Q):Q(P,Q))  \\
    &=&  (H^d P(p,q):H^d Q(p,q)) \\
    &=&  (H^d H(p,q) p(p,q): H^d H(p,q) q(p,q)) \\
    &=&  H^d \phi^*(H) \phi\circ\phi.
\end{eqnarray*}
Since the map $\phi$ is non-degenerate, we will never have
$\phi\circ\phi(z:w) = (0:0)$.  However, we have 
$H(p,q) \equiv 0$ if and only if $\phi(z:w)
=(\alpha:\beta)\in\P^1$ for all $(z:w)\in\P^1$ and 
$H(\alpha,\beta)=0$.
This exactly characterizes the set $I(d)$.  
Thus, the above gives the formula of the second iterate
for $f\not\in I(d)$, and the second iterate is undefined for 
$f\in I(d)$. 

An easy inductive argument 
gives the general form of the iterate $f^n$ for all $f\not\in I(d)$. 
Since the formula for $f^n$  does not vanish identically
for any $f\not\in I(d)$,  
the indeterminacy locus of $\Phi_n$ must be 
contained in $I(d)$ for each $n\geq 3$.  However, by Lemma
\ref{descending}, the chain of ideals defined by the iterate
maps is descending.  Thus, if the coordinate functions of 
$\Phi_2$ vanish simultaneously along $I(d)$, the coordinate 
functions of $\Phi_n$ vanish simultaneously along $I(d)$
for all $n\geq 2$.  Therefore, the indeterminacy locus is 
$I(d)$ for all $n$. 
\qed

\medskip
Note that indeterminacy of the iterate map along $I(d)$ implies
that $\Phi_n$ can not be extended continuously from $\Rat_d$ 
to any point 
$g\in I(d)$.  See the two examples in Section \ref{Thm2}.
Observe also that Lemma \ref{iterate} is a statement about the 
indeterminacy locus as a {\em set}.  Scheme-theoretically, the 
indeterminacy depends on the iterate $n$.  

As holomorphic correspondences in $\P^1\times\P^1$, 
the elements of $I(d)$ each consist of a flat horizontal 
component and a collection of vertical components, one 
of which must intersect the horizontal graph 
on the diagonal.   The second iterate 
of a point $f\in\Ratbar_d$ as a correspondence is given by
  $$\Gamma_f\circ\Gamma_f = \{(z_1,z_2)\in\P^1\times\P^1: 
      z_2 = f(y) \mbox{ and } y=f(z_1) \mbox{ for some } 
       y\in\P^1\}.$$
With this notion of iteration, it is easy to see that the 
indeterminacy locus $I(d)$ satisfies
  $$I(d) = \{f\in \P\H^0(\P^1\times\P^1,\O(d,1)):
             \Gamma_f\circ\Gamma_f = \P^1\times\P^1\}.$$ 
For $f\not\in I(d)$, the composition $\Gamma_f\circ\Gamma_f$
is the zero-locus of a well-defined section 
of the line-bundle $\O(d^2,1)$ over $\P^1\times\P^1$, 
and Lemma \ref{iterate}
describes its graph and the location of the vertical components.

Note, for example, if $f = H\phi$ is degenerate of degree $d$, 
then the holes of the third
iterate of $f$ are the holes of $f$ at depth $d^2$, the preimages
by $\phi$ of these holes at depth $d$, and finally the preimages
by $\phi^2$ of the holes of $f$.  Comparing the 
iterate formula  
with the definition of the measure $\mu_f$ 
given in Section \ref{definitions} (and Lemma \ref{pointmass}), 
we obtain the following immediate corollary.  

\begin{cor} \label{holemass}
Let $f\not\in I(d)$ be degenerate.  
If $d_z(f^n)$ denotes the depth of $z$ as a hole of $f^n$, then 
the $d_z(f^n)/d^n$ forms a non-decreasing sequence in $n$, and 
  $$\mu_f(\{z\}) = \lim_{n\to\infty} \frac{d_z(f^n)}{d^n}.$$
Furthermore, $\mu_f = \mu_{f^n}$ for all $n\geq 1$.  
\end{cor}

%\begin{lemma} \label{graph}
%For each $n\geq 2$, the closure of the graph of $\Phi_n$ in 
%$\P^{2d+1}\times\P^{2d^n+1}$ is equal to the closure of the 
%graph of $\Phi_n|\Rat_d$. 
%\end{lemma}
%
%\proof
%Let $\Gamma$ be the closure of the graph of $\Phi_n$; that is,
%  $$\Gamma = \overline{\{ (f,\Phi_n(f)): f\in\P^{2d+1}-I(d) \}} 
%         \subset \P^{2d+1}\times\P^{2d^n+1}. $$
%Then $\Gamma$ is irreducible (since $\P^{2d+1}$ is irreducible).
%Therefore, $\Gamma$ is equal to the closure of the graph of 
%$\Phi_n$ restricted to any Zariski open subset of $\P^{2d+1}-I(d)$.
%See, for example, \cite[Ch.? Thm?]{Whitney:varieties}.
%\qed

%%%%%%%%%%

\section{The dynamics of a degenerate map}  \label{dynamics}

In this section we define the Julia set of a degenerate map 
$f\not\in I(d)$ and 
relate it to the support of the measure $\mu_f$.  We explain 
how $\mu_f$ is the weak limit of pull-backs by $f$ of any probability
measure on the Riemann sphere.  Also, it is possible to define
the escape rate function in $\C^2$ for every degenerate 
$f\not\in I(d)$, and we will see that it is a potential for the atomic
measure $\mu_f$.

Let $f = H\phi \not\in I(d)$ be degenerate.  As for rational maps, we define
the Fatou set $\Omega(f)$ as the largest open set on which the 
iterates of $f$ form a normal family.  
Care must be taken in this definition since we require, first, that the
iterates $f^n$ be well-defined for each $n$.  Thus, the 
family $\{f^n|U\}_n$ can not be normal 
if $h\in f^n(U)$ for some $n\geq 1$ and some hole $h$ of $f$.
With this definition, we let $J(f)$ be the complement of $\Omega(f)$.

Let us assume for the moment that $\deg \phi >0$.  
By the definition of the Julia set $J(f)$, it is clear that 
 $$J(f) = J(\phi) \cup  \overline{ \bigcup_{n=0}^\infty 
\bigcup_i \phi^{-n}(h_i) },$$
where $\{h_i\}\subset\P^1$ is the set of holes of $f$.
Recall that a point $z\in \P^1$ is said to be exceptional for $\phi$ if 
its grand orbit is finite.  
If at least one of the holes $h_i$ is non-exceptional
for the map $\phi$, then $J(\phi)$ is contained in the 
closure of the union of the preimages of $h_i$.  Examining again the 
definition of the measure $\mu_f$, we see that its support must
be the closure of the union of all preimages of the holes of $f$.
When $\deg \phi =0$, it makes sense to set $J(\phi)= \emptyset$.
We have proved the following.  

\begin{prop}
Let $f = H\phi \not\in I(d)$ be degenerate.  If one of the holes
of $f$ is non-exceptional for $\phi$, then 
$\supp \mu_f = J(f)$. 
If each hole of $f$ is exceptional for $\phi$, then $\supp \mu_f$ is
contained in the exceptional set of $\phi$.
\end{prop}

\noindent
Note that even if $J(\phi)\subset\supp\mu_f$, it can happen
that $\mu_f(J(\phi))=0$, as all holes of $f$ may lie in the Fatou
set of $\phi$.  

\bigskip\noindent
{\bf Pullbacks of measures by degenerate maps}.  
The holes of a degenerate map $f$ are identified with the 
vertical components of the holomorphic correspondence of $f$,
as described in Section \ref{definitions}.  
The degenerate $f$ should be interpreted as sending each of 
its holes over the whole of $\P^1$ with appropriate multiplicity (the 
depth of the hole).  In particular, pull-backs
of measures can be appropriately defined, at least when $\deg
\phi\not=0$, by integration over the fibers, as
 $$\< f^*\mu,\psi\>  
= \< \phi^*\mu,\psi\> + \sum_i \< \delta_{h_i},\psi\>
= \int \sum_{\phi(y)=z}\psi(y) ~d\mu(z) 
+ \sum_i \psi(h_i),$$
where $\{h_i\}\subset\P^1$ is the set of holes of $f$ and 
$\psi$ is any continuous function on $\P^1$.   All sums are
counted with multiplicity.  When $\deg\phi=0$ but 
$f\not\in I(d)$, we can set 
 $$\< f^*\mu,\psi \> = \sum_i \psi(h_i).$$  

The following proposition implies that for each $f\in\del\Rat_d - I(d)$, 
there is a unique fixed point of the
operator $\mu\mapsto f^*\mu/d$ in the space of probability measures. 
In particular, a degenerate map not in $I(d)$ has no exceptional points.

\begin{prop}
For any degenerate $f=\not\in I(d)$ and any probability measure $\mu$ on 
$\P^1$, we have $f^{n*}\mu/d^n \to \mu_f$
weakly as $n\to\infty$.
\end{prop}

\proof
Write $f = H\phi$.  
As the degree of $\phi$ is strictly less than $d$, we have
  $$ \frac{1}{d^n} \left| \< \phi^{n*}\mu, \psi \> \right| \to 0,$$
for all test functions $\psi$.  From Corollary \ref{holemass}, 
the depths of  the holes of the iterates of $f$ limit on the mass 
$\mu_f$.  
\qed

\bigskip\noindent{\bf Escape rate functions.}
The escape rate function of a rational map $f\in\Rat_d$, $d\geq 2$, 
is defined by 
\begin{equation} \label{G_F}
 G_F(z,w) = \lim_{n\to\infty} \frac{1}{d^n} \log \|F^n(z,w)\|,
\end{equation}
where $F:\C^2\to\C^2$ is a homogeneous polynomial map such 
that $\pi\circ F = f\circ\pi$.  Here, $\pi:\C^2-0\to\P^1$ is the 
canonical projection and $\|\cdot\|$ is any norm
on $\C^2$.  If $F_1$ and $F_2$ are two lifts of $f$ to 
$\C^2$ (so that necessarily $F_2=\alpha F_1$ for some $\alpha\in\C^*$), 
then $G_{F_1}-G_{F_2}$ is constant. 
The escape rate function is a potential for the 
measure $\mu_f$ in the sense that $\pi^*\mu_f = dd^c G_F$  
\cite[Thm 4.1]{Hubbard:Papadopol}.  We use the notation
$d = \del + \delbar$ and $d^c = i(\delbar-\del)/2\pi$.  

The proofs of the following proposition and Corollary
\ref{limit G_F} rely on the isomorphism 
between the space of probability measures on $\P^1$ and (normalized) 
plurisubharmonic functions $U$ on $\C^2$ such that 
$U(\alpha z) = U(z) + \log |\alpha|$ for all $\alpha\in\C^*$ 
\cite[Thm 5.9]{Fornaess:Sibony}.  
The isomorphism is given by $\mu= \pi_* dd^cU$, 
from potential functions to measures.

\begin{prop} \label{escape}
For $d\geq 2$, the escape rate function $G_F$ exists for 
each degenerate $f\not\in I(d)$
and satisfies $\pi^*\mu_f = dd^c G_F$.
\end{prop}

We first need a lemma on M\"obius transformations.  Let 
$\sigma$ denote the spherical metric on $\P^1$ and 
$\dist_\sigma$ the associated distance function.  

\begin{lemma} \label{Mobius}
Let $E\subset\P^1$ be a finite set and $E(r) = \{z\in\P^1: 
\dist_\sigma(z,E) \leq r\}$.  For each M\"obius transformation 
$M\in\Rat_1$, there exists $r_0>0$ such that 
  $$\bigcup_{k\geq 0} M^k(E(r^k)) \not= \P^1$$
for all $r<r_0$.  
\end{lemma}

\proof
By choosing coordinates on $\P^1\iso\Chat$, 
we can assume that 
$M$ has the form $M(z) = z+1$ or $M(z) = \lambda z$ for 
$\lambda\in\overline{\D}$.  In the new coordinate system, 
the spherical metric is comparable to the pull-back of the 
given metric $\sigma$ by the coordinate map.  

When $|\lambda|=1$, the statement is obvious for $r_0$ 
sufficiently small.  When $|\lambda|<1$, we need only consider
the case when $\infty\in E$.  For $r$ small, a spherical disk 
of radius $r$ around $\infty$ is comparable in size to the 
complement of the Euclidean disk of radius $1/r$ centered at 0, 
so we need to choose $r_0 < |\lambda|$.

Finally, suppose $M(z) = z+1$.  Again, we need only consider 
the case when $\infty\in E$.  It is clear that the point $M^k(0)=k$,
for example, remains inside the Euclidean disk of radius 
$r^{-k}$ for all $k$ if $r<1$.  Thus, $M^k(0)$ is outside the 
spherical disk of radius $r^k$ about $\infty$ for all $k$. 
We can therefore choose any $r_0<1$.  
\qed

\medskip\noindent{\bf Proof of Proposition \ref{escape}}.
Let $f = H \phi\not\in I(d)$ be degenerate.  
Expressing $f$ in homogeneous coordinates 
defines a polynomial map $F:\C^2\to\C^2$ of (algebraic) degree 
$d$ such that $\pi\circ F=f\circ\pi$ where defined.  In particular,
$F$ vanishes identically along the lines $\pi^{-1}(h)$ for 
each hole $h$ of $f$.  We aim to define $G_F$ 
by equation (\ref{G_F}) as for non-degenerate 
maps, so we need to show that the limit exists.  

Write $F = (P,Q) = H \Phi$, where $H= \gcd(P,Q)$ and $\Phi
= (P/H, Q/H)$
is a non-degenerate homogeneous polynomial map of degree $e
< d$ (so that $\Phi^{-1}(\{0\})=\{0\}$).  In fact, $\Phi$ is the
the map $\phi$ expressed in homogeneous coordinates.    

The iterate formula for 
$f$ in Lemma \ref{iterate} holds also for $F$ so that 
\begin{equation} \label{G_n}
 G_n(x) := \frac{1}{d^n}\log\|F^n(x)\| 
      = \sum_{k=0}^{n-1} \frac{1}{d^{k+1}} \log |H(\Phi^k(x))|
           + \frac{1}{d^n}\log\|\Phi^n(x)\|
\end{equation}
for all $x\in\C^2$.  

Suppose first that $e = \deg \phi = 0$, so that $F(z,w) 
= (aH(z,w),bH(z,w))$ and $H(a,b)\not=0$. The above expression 
for $G_n$ reduces to  
 $$G_n(z,w) = \frac{1}{d} \log|H(z,w)| + 
              \sum_1^{n-1} \frac{1}{d^{k+1}} \log|H(a,b)| 
              + \frac{1}{d^n} \log\|(a,b)\|,$$
which converges to 
 $$G_F(z,w) = \frac{1}{d} \log|H(z,w)| + 
              \frac{1}{d(d-1)} \log|H(a,b)|$$
locally uniformly on $\C^2-0$ as $n\to\infty$.  
Furthermore, this function $G_F$
is clearly a potential for the atomic measure 
  $$\mu_f = \frac{1}{d} \sum_{z\in\P^1:H|\pi^{-1}(z)=0} \delta_z$$
on $\P^1$, where the zeros of $H$ are counted with multiplicity.

Now suppose that $e = \deg\phi>0$.  Then there exists
a constant $K>1$ so that for all $x\in\C^2$, 
  $$K^{-1}\|x\|^e \leq \|\Phi(x)\| 
               \leq K \|x\|^e$$
and therefore if $x\not=0$, 
  $$|\log\|\Phi(x)\|| \leq e |\log\|x\|| + \log K.$$
Replacing $x$ with the iterate $\Phi^{n-1}(x)$ we obtain 
by induction on $n$, 
  $$|\log\|\Phi^n(x)\|| \leq e^n |\log\|x\|| 
                 + (1+e+\cdots + e^{n-1})\log K.$$
Dividing by $d^n$ gives 
\begin{equation} \label{to 0}
 \frac{1}{d^n} \log\|\Phi^n(x)\| \to 0
\end{equation}
locally uniformly on $\C^2-0$ as $n\to\infty$,
since $e<d$.  Similarly, the quantity
\begin{equation} \label{bound}
 \sum_0^{n-1} \frac{1}{d^{k+1}} \log \|\Phi^k(x)\| 
      - \sum_0^{n-1} \frac{e^k}{d^{k+1}} \log \|x\|
\end{equation}
is uniformly bounded in $n$ on $\C^2-0$.

Consider the plurisubharmonic function 
  $$g_n(x) = 
   \sum_{k=0}^{n-1} \frac{1}{d^{k+1}} \log |H(\Phi^k(x))|,$$
on $\C^2$.  Notice
that $g_n$ defines a potential function for the atomic measure 
  $$\mu_n = \sum_i \sum_{k=0}^{n-1} \frac{1}{d^{k+1}} 
                \sum_{z:\phi^k(z)=h_i} \delta_z,$$
where the $h_i$ are the holes of $f$, counted with multiplicity;
that is, $\pi^* \mu_n = dd^c g_n$ on $\C^2-0$.  
Note also that the measure $\mu_n$ has total mass $1-(e/d)^n$,
and $g_n$ scales by 
  $$g_n(\alpha x) = g_n(x) + \left( 1-(e/d)^n\right)
                      \log|\alpha|$$
for all $\alpha\in\C^*$.  The measures $\mu_n$ converge weakly 
to $\mu_f$ in $\P^1$.  By (\ref{G_n}) and (\ref{to 0}), for any 
$\eps>0$, we have 
\begin{equation}  \label{G_n-g_n}
  |G_n(x) - g_n(x)| < \eps
\end{equation}
for all sufficiently large $n$, locally uniformly in $\C^2-0$.
We will show that the functions $g_n$ converge in $L^1_{loc}$ to
the unique potential function of $\mu_f$.  

If the sequence $g_n$ is uniformly bounded above on compact 
sets and does not converge to $-\infty$ locally uniformly, then
some subsequence converges in $L^1_{loc}$ 
\cite[Thm 5.1]{Fornaess:Sibony}.  For an upper bound, note first
that
\begin{equation} \label{sup H}
  \sup\{\log|H(x)|: \|x\|\leq 1\} < \infty.  
\end{equation}
If for each $x\not=0$ in $\C^2$, we set $x^1 := x/\|x\|$, then 
$H(x) = \|x\|^{d-e} H(x^1)$, and therefore, 
\begin{equation} \label {g_n}
  g_n(x) = \sum_0^{n-1} \frac{d-e}{d^{k+1}} 
     \log\|\Phi^k(x)\| 
    + \sum_0^{n-1} \frac{1}{d^{k+1}} \log |H(\Phi^k(x)^1)|.
\end{equation}
The bound on (\ref{bound}) together with (\ref{sup H}) show that 
$\{g_n\}$ is uniformly bounded above on compact sets.  

To obtain a convergent subsequence of $\{g_n\}$ in $L^1_{loc}$, 
it suffices now (by \cite[Thm 5.1]{Fornaess:Sibony})
to show that $g_n\not\to -\infty$ uniformly on compact sets.
If  $g_{n_j}$ converges to $v$ in $L^1_{loc}$,
then by \cite[Thm 5.9]{Fornaess:Sibony}, $v$ is the potential 
function for the measure $\mu_f$,
unique up to an additive constant.  To conclude, therefore, 
that the full sequence $g_n$ converges so that $\pi^*\mu_f
= dd^c G_F$, it suffices to show that there exists a single
point $x\in\C^2-0$ for which $\lim_{n\to\infty} g_n(x)$ 
exists and is finite.  

For $e\geq 2$, choose any point $x\in\C^2-0$ which 
is periodic for $\Phi$ and whose orbit does not intersect 
the complex lines $\pi^{-1}(h_i)$ over the holes of $f$.  Then the
orbit $\Phi^k(x)$, $k\geq 0$, remains a bounded distance away from 
the zeroes of $H$.  Therefore, 
$\log|H(\Phi^k(x))|$ is bounded above and below, so that 
the definition of $g_n$ together with (\ref{G_n-g_n}) imply that 
  $$G_F(x) = \lim_{n\to\infty} G_n(x) = \lim_{n\to\infty} g_n(x)$$
exists and is finite.
Therefore, $G_F\in L^1_{loc}$ is the potential function for
$\mu_f$.  

For $e=1$ a further estimate is required.  The map 
$\phi$ on $\P^1$ is a Mobius transformation.  Let $E$ be 
the set of zeroes of $H$ projected to $\P^1$.  Let $\sigma$
denote the spherical metric on $\P^1$, and observe that 
there exists $C>1$ and a positive integer $m$ such that 
\begin{equation} \label{H(x)}
  C \geq |H(x)| \geq C^{-1} \dist_\sigma(\pi(x),E)^m 
    \mbox{ for all } \|x\|=1.
\end{equation}
By Lemma \ref{Mobius} applied to $M=\phi^{-1}$, 
the set $A = \P^1 - \cup_{k\geq 0} \phi^{-k} E(r^k)$ 
is non-empty for some $r>0$.
Fix $z\in A$ and choose $Z\in\pi^{-1}(z)$ with $\|Z\|=1$.  
Using (\ref{H(x)}), our choice of $Z$ implies that
  $$\log C \geq \log|H(\Phi^k(Z)^1)| \geq mk\log r - \log C,$$
for all $k\geq 1$, where $\Phi^k(Z)^1 = \Phi^k(Z)/\|\Phi^k(Z)\|$.
Note that since $\Phi$ is linear, there is a 
constant $C'>1$ such that 
  $$|\log\|\Phi^k(Z)\||\leq k \log C'.$$  
Combining these estimates, we obtain
$$ \sum_n^{m} \left| \frac{d-1}{d^{k+1}} \log\|\Phi^k(Z)\|  
               + \frac{1}{d^{k+1}} \log 
                      |H(\Phi^k(Z)^1)| \right|
   \leq C'' \sum_n^m \frac{k}{d^k},  $$
for some constant $C''$ 
which implies, by (\ref{g_n}), that $g_n(Z)$ has a finite 
limit as $n\to\infty$.   
\qed

\medskip
As a corollary to Theorem 1(a), we obtain

\begin{cor} \label{limit G_F}
Suppose that a sequence $f_k$ in $\Rat_d$ converges in $\Ratbar_d$
to $f\not\in I(d)$.  Then for 
suitably normalized lifts $F_k$ of $f$ to $\C^2$ and $F$ of $f$, 
the escape rate functions $G_{F_k}$ converge to $G_F$ in $L^1_{loc}$.
\end{cor}

\proof
Again by \cite[Thm 5.9]{Fornaess:Sibony}, weak convergence
of measures $\mu_{f_k}\to\mu_f$  by Theorem 1(a) implies
that the potentials converge 
in $L^1_{loc}$.  A normalization is required 
to guarantee convergence; it suffices to choose the unique lifts 
$F_k$ of $f_k$ such that $\sup\{G_{F_k}(x): \|x\|=1\} = 0$.  
\qed

%%%%%%%%

\section{Proof of Theorem 1}  \label{Thm1}

Here we provide the statements needed for the proof of
Theorem 1.  The argument relies on a fundamental fact about holomorphic 
functions in $\C$:  a proper holomorphic function from a domain $U$ 
to a domain $V$ in $\C$ has a well-defined degree.  This together 
with the invariance property of the maximal measure, $f^*\mu_f/d = 
\mu_f$ for all $f\in\Rat_d$, and the continuity of the 
iterate map away from $I(d)$ will give Theorem 1.  
For simplicity, we will regularly identify the point $(z:w)\in\P^1$ with 
$z/w\in \overline{\C}$.  

\begin{lemma} \label{0s and poles}
Suppose a sequence $f_k\in\Rat_d$ converges
in $\Ratbar_d$ to degenerate $f=(P:Q)$ 
and $f$ has a hole at $h$ of depth $d_h$.  
If neither $P$ nor $Q$ are $\equiv 0$, then any 
neighborhood $N$ of $h$ contains at least $d_h$ 
zeroes and poles of $f_k$ (counted with multiplicity) for all 
sufficiently large $k$.  
\end{lemma}

\proof
As the coefficients of $f_k=(P_k:Q_k)$ converge to those of $f$, then 
so must the 
roots of the polynomials $P_k$ and $Q_k$ converge to those of $P$ and 
$Q$.  If the degenerate map $f$ has a hole at $h\in\P^1$ of depth $d_h$, 
then at least $d_h$ roots of $P$ and $d_h$ roots of $Q$ must limit on 
$h$.
\qed

\begin{lemma} \label{convergence}
Suppose $f_k\in\Rat_d$ converge to $f=H\phi$ in $\Ratbar_d$.  Then 
as maps, $f_k\to\phi$ locally uniformly on $\P^1-\{h_i\}_i$ 
where the $h_i$ are the holes of $f$.  
\end{lemma}

\proof
Write $f_k=(P_k:Q_k)$ and $f=(P:Q)$.  By changing coordinates, we may 
assume that neither $P$ nor $Q$ are $\equiv 0$ and also that no holes 
lie at the point $\infty = (1:0)$.  Writing $f(z:w) = 
(H(z,w)p(z,w):H(z,w)q(z,w))$, where $H = \gcd(P,Q)$ we may assume 
that $H$ is monic as a polynomial in $z$ of degree $d-e$.  
Fix an open set $U\subset\P^1$ 
containing all holes of $f$.  By the previous lemma there are 
homogeneous factors 
$A_k(z,w)$ of $P_k$ and $B_k(z,w)$ of $Q_k$ of degree $d-e$ with 
all roots inside $U$.  As polynomials of $z$, we may assume that $A_k$ 
and $B_k$ are monic, and thus $A_k$ and $B_k$ both tend to $H$ as 
$k\to\infty$.  On the compact set $\P^1-U$, the ratio $A_k/B_k$ 
must tend uniformly to 1 as $k\to\infty$, so we have 
$f_k=(P_k:Q_k)$ limiting on $\phi=(p:q)$ 
uniformly on $\P^1-U$.
\qed

\medskip
For the proof of the theorem, we will need a uniform version of Lemma
\ref{0s and poles}; namely, there should be preimages of almost every 
point (with respect to the maximal measure of $f_k$) inside a small 
neighborhood of the holes.  As we shall see, this can be done when the 
limit map is not in $I(d)$.  Uniformity fails in general, and this
failure leads to the discontinuity of the measure in Theorem 2.

\begin{prop} \label{generic case}
Suppose the sequence $f_k$ in $\Rat_d$ converges to degenerate 
$f\not\in I(d)$.
If $f$ has a hole at $h$ of depth $d_h$, then any weak limit $\nu$ of 
the maximal
measures $\mu_{f_k}$ of $f_k$ must have $\nu(\{h\}) \geq d_h/d$.  
\end{prop}

\begin{prop} \label{special case}
Suppose the sequence $f_k$ in $\Rat_d$ converges to 
$f = H\phi \in  I(d)$, 
If $f$ has a hole of depth $d_c$ at $c$ where $\phi\equiv c$,
then any weak limit $\nu$ of the maximal measures 
$\mu_{f_k}$ of $f_k$ must 
have $\nu(\{c\}) \geq d_c/(d_c+d)$.  
\end{prop}

The proof of Proposition \ref{generic case} follows from the following two 
lemmas.

\begin{lemma} \label{A}
Suppose under the hypotheses of Proposition \ref{generic case} that 
$\phi$ is non-constant.  Then for any neighborhood $N$ of the hole 
$h$, there exists an $M>0$ such that 
  $$  \#\{f_k^{-1}(a)\cap N\} \geq d_h,$$
for all $a\in\P^1$ and all $k>M$, where the preimages are counted with 
multiplicity.  
\end{lemma}

\proof
Suppose that in local coordinates at $h$ and $\phi(h)$, we can 
write $\phi(z) = c z^m + O(z^{m+1})$, $m>0$.   Choose a disk $D$
around $\phi(h)$ small enough that 
\begin{itemize}
\item[(i)]  $D$ does not contain both 0 and $\infty$,
\item[(ii)]  the component of $\phi^{-1}(D)$ containing $h$ 
is a disk inside $N$, and
\item[(iii)]  the component of $\phi^{-1}(D)$ containing $h$ maps 
               $m$-to-1 over $D$.
\end{itemize}
Let this component of the preimage of $D$ be denoted by $E$.

By uniform convergence of $f_k$ to $\phi$
away from the holes of $f$ (Lemma \ref{convergence}), 
for all sufficiently large $k$, $f_k$ maps a curve
close to $\del E$ $m$-to-1 over $\del D$, and by Lemma \ref{0s and 
poles}, $d_h$ zeros or poles lie very close to $h$.  (Note that the 
hypothesis of Lemma \ref{0s and poles} is automatically satisfied when 
$\phi$ is non-constant.)  Let $E_k$ denote
the disk containing $h$ bounded by the component of $f_k^{-1}(\del D)$
which is very close to $\del E$.
Consider the preimage $A_k = f_k^{-1}(\P^1-D) \cap E_k$.  As $f_k$ is
proper on $A_k$, it has a well-defined degree.  Counting zeros or poles
in $A_k$, we find the degree is $d_h$.  The map $f_k$ then has degree
$d_h$ also on the boundary of $A_k$.  Now, $f_k$ is also proper on the 
complement of $A_k$ in $E_k$.  Counting degree on its boundary,
$d_h+m$, we find that $f_k$ has at least $d_h$ preimages of {\em all} 
points of the sphere inside $N$.   \qed

\begin{lemma}  \label{B}
Suppose under the hypotheses of Proposition \ref{generic case} that 
$\phi$ is constant.  Then for any neighborhood $N$ of $h$, there 
exists an $M>0$ so that 
  $$\# \{f_k^{-1}(a)\cap N\} \geq d_h,$$
for all $k>M$ and all $a\in \supp(\mu_{f_k})$, where the preimages are 
counted with multiplicity.  
\end{lemma}

\proof
Suppose that $\phi\equiv c$.  By changing coordinates if necessary,
we can assume that the point $c$ is neither $0=(0:1)$ nor $\infty=(1:0)$
so that the hypothesis of Lemma \ref{0s and poles} is satisfied.  By
assumption, $c$ is not one of the holes of $f$.  Let $D$ be a disk
around $c$ so that all holes and one of 0 or $\infty$ lie outside $D$.  
Let $B$ be a ball around $h$ contained in $N$.  For all large $k$, $f_k$
has $d_h$ zeros and poles inside $B$ and $f_k$ maps the complement of 
$B$ (minus a neighborhood of other holes) to $D$.

It is clear that for these large $k$, $D$ does not interesect the
Julia set of $f_k$, since $f_k(D)$ is contained in $D$ so the
iterates must form a normal family on $D$.  

Consider the preimage $A_k = f_k^{-1}(\P^1-D) \cap B$.  The map
$f_k$ is proper on $A_k$ and has a well-defined degree.  Counting
zeros or poles, this degree is $d_h$.  Since $\supp\mu_{f_k}$ lies
in $\P^1-D$, the lemma is proved.  \qed

\medskip\noindent
{\bf Proof of Proposition \ref{generic case}}.  
This follows immediately from the invariance property of 
$\mu_{f_k}$:  
   $$f_k^* \mu_{f_k} = d \, \mu_{f_k}.$$
Fix a disk $D$ around the hole $h$.  Choose a bump function $b$ which 
is 1 on a disk of half the radius and 0 outside $D$.  Then 
 $$\mu_{f_k}(D) \geq \int b \,d\mu_{f_k} 
 = \frac{1}{d} \int \sum_{f_k(x)=y} b(x) \,d\mu_{f_k}(y).$$
By Lemmas \ref{A} and \ref{B}, for all sufficiently large $k$, the sum in 
the integrand is $\geq d_h$ for every $y$ in the support of $\mu_{f_k}$. 
Taking limits and letting $D$ shrink down to $h$ gives the result.  
\qed

\medskip\noindent
{\bf Proof of Proposition \ref{special case}}.
By changing coordinates if necessary, we can assume that the constant 
$c$ is neither $0=(0:1)$ nor $\infty=(1:0)$ so that the hypothesis 
of Lemma \ref{0s and poles} is satisfied.  
Let $D$ be a disk around $c$ which does not contain both 
0 and $\infty$.  Let $N$ be a neighborhood of the 
holes of $f$ such that $\overline{N}\cap\del D = \emptyset$. 
Now choose $M$ large enough that 
$f_k(\P^1-N)\subset D$ for all $k>M$.  Let $A_k = 
f_k^{-1}(\P^1-D) \cap D$.  The map $f_k$ is proper on $A_k$ 
and has at least $d_c$ zeroes or poles, so it is at least $d_c$-to-1 
over the complement of $D$.  

Let $\nu$ be any subsequential weak limit of $\mu_{f_k}$.  
Let $\nu_D = \nu(\P^1-\overline{D})$.  If $\nu_D=0$ and this holds for 
all $D$, then $\nu(\{c\})=1$ and the proposition is proved. 
Generally, for any $\eps>0$ we have 
$\mu_{f_k}(\P^1-\overline{D})\geq \nu_D-\eps$
for all large $k$.  By the invariance property of $\mu_{f_k}$ (as in 
the proof of the previous proposition), $\mu_{f_k}(D)\geq 
(\nu_D-\eps)d_c/d$.
Thus $\nu(D)\geq (\nu_D-\eps)d_c/d$.  Since $\eps$ is arbitrary and 
since $\nu_D$ increases (to some value $\nu_0\leq 1$) as $D$ shrinks,
we obtain $\nu(\{c\}) \geq \nu_0d_c/d$.  

On the other hand, 
 $$1 = \nu_0 + \nu(\{c\}) \geq \nu_0 + \nu_0d_c/d = \nu_0(1 + d_c/d),$$
and therefore $\nu_0 \leq d/(d_c+d)$.  Consequently, $\nu(\{c\}) \geq
d_c/(d_c+d)$.  
\qed

\bigskip\noindent
{\bf Proof of Theorem 1}.  
By Lemma \ref{iterate}, the iterate map $\Phi_n$ is continuous on 
$\Ratbar_d-I(d)$ for every $n$.  
Thus, if $f_k\to f\not\in I(d)$, then $f_k^n\to f^n$ 
where the iterate of a degenerate map is described explicitly in 
Lemma \ref{iterate}.  From Corollary \ref{holemass}, 
we know that $1/d^n$ multiplied by the depth of a hole of 
$f^n$ can only increase as $n\to\infty$.  Also,
the depths limit on the mass given by $\mu_f$.  

Since the maximal measure for a rational map is the same as the measure 
for any iterate, Proposition \ref{generic case}
implies that any subsequential limit $\nu$ of the measures has {\em
at least} the correct mass on all the points in $\supp{\mu_f}$.  On 
the other hand, these masses sum to 1, and the measure is a 
probability measure, so in fact, $\nu = \mu_f$.  This proves part (a).
Part (b) is exactly the statement of Proposition \ref{special case}. 
\qed

%%%%%%%%%%

\section{Examples and proof of Theorem 2} \label{Thm2}

In this section we complete the proof of Theorem 2.  We begin with 
some examples demonstrating the discontinuity 
of the iterate maps and, consequently, the discontinuity of the 
map of measures $f\mapsto\mu_f$ at each point in $I(d)$. 
Example 1 realizes the lower bound of Theorem 1(b) when the
depth $d_c$ is 1.  
A point $(z:w)\in\P^1$ will be regularly identified
with the ratio $z/w$ in $\overline{\C}$.  

\medskip\noindent
{\bf Example 1}.  Let $g = (w P(z,w):0)\in I(d)$ where $P$ is homogeneous
of degree $d-1$, $P(0,1)\not=0$, $P(1,0)\not=0$, and $P$ is 
monic as a polynomial in $z$.  Then $g$ has a hole of depth 1 at
$\infty$ and no holes at 0.  For each $a\in\C$ and 
$t\in\D^*$, consider 
 $$ g_{a,t}(z:w) := (a t z^d + w P(z,w): t z^d) \in \Rat_d,$$
so that $g_{a,t} \to g$ in $\Ratbar_d$ as $t\to 0$.  The maps  
$g_{a,t}$ all have a critical point at $z=0$ of multiplicity $d-1$, and
the other critical points are at the $d-1$ solutions to 
$zP'(z,1) - dP(z,1) = 0$, independent of both $a$ and $t$.  
For each $a\in\C$, $g_{a,t}$ converges to the constant 
$\infty$ as $t\to 0$, uniformly 
away from $\infty$ and the roots of $P$ (by Lemma \ref{convergence}).  
The second iterate $\Phi_2(g_{a,t})$ has the form
 $$(aw^dP(z,w)^dt + z^dw^{d-1}P(z,w)^{d-1}t + O(t^2) :
                      w^dP(z,w)^dt + O(t^2)),$$
and taking a limit as $t\to 0$, we obtain 
 $$\Phi_2(g_{a,t}) \to f_a := (w^{d-1}P(z,w)^{d-1}(awP(z,w)
                    + z^d) : w^dP(z,w)^d). $$
Thus the second iterates converge (uniformly, away from the holes of 
$f_a$) to the map $\phi_a$ given (in coordinates on $\C$) by
  $$\phi_a(z) = \frac{aP(z) + z^d}{P(z)}.$$
Recall that we are assuming $P(0)\not=0$ so that $\phi_a$ 
is a non-degenerate rational map of degree $d$ for all 
$a\in\C$.  As $P$ is monic and of degree $d-1$,
it is clear that each $\phi_a$ has a parabolic fixed point 
at $\infty$.
We see immediately that the limit as $t\to 0$
of second iterates of $g_{a,t}$ depends 
on the direction of approach, parameterized here by $a\in\C$.
The degenerate limits $f_a$ all have holes at $\infty$ of depth $d-1$
and holes at the roots of $P$ each of depth $d-1$.

The degenerate maps $f_a\in\Ratbar_{d^2}$ do not 
lie in $I(d^2)$, so we are able to compute the limiting measures 
of $\mu_{g_{a,t}}$ as $t\to 0$ from Theorem 1(a).  
Of course, the maximal measure for 
$g_{a,t}$ coincides with that of its second iterate, so
$\mu_{g_{a,t}} \to \mu_{f_a}$ weakly for each $a\in\C$ as 
$t\to 0$.  

The measures $\mu_{f_a}$ cannot be the same for all $a\in\C$:  
the holes are the same for each $a$ but 
the preimages of the roots of $P$ by $\phi_a$
depend on $a$, and these are atoms of $\mu_{f_a}$.    For 
example, suppose $\alpha\not=0$ is a simple root of $P(z)$.  Then 
for $a=\alpha$, the 
$d$ solutions to $\phi_\alpha(z) = \alpha$ are all at 0, so that by Lemma 
\ref{pointmass},
 $$\mu_{f_\alpha}(\{0\}) = \frac{1}{d^2} \sum_{n=1}^\infty 
                \frac{d(d-1)}{d^{2n}} = \frac{1}{d(d+1)}.$$
On the other hand, for the generic $a\in\C$, the $\phi_a$-orbit 
of the point 0 never intersects the roots of $P$, so that 
$\mu_{f_a}(\{0\}) =0$.

Lastly, since the limiting maps $f_a$ have holes of depth
$d-1$ at $\infty$ and $\phi_a(\infty)=\infty$ for each 
$a\in\C$, we can compute
easily from Lemma \ref{pointmass} that 
 $$\mu_{f_a}(\{\infty\}) = \frac{1}{d^2} 
            \sum_0^\infty \frac{d-1}{d^{2n}}  = \frac{1}{d+1},$$
for all $a\in\C$.
As the degeneration of $g_{a,t}$ to $g$ develops a hole of depth
$d_{\infty} = 1$ at $\infty$, we see that this family achieves
the lower bound of Theorem 1(b).

\medskip\noindent{\bf Example 2}.  
Let $g = (w^kP(z,w):0)\in I(d)$ where $P$ is homogeneous of 
degree $d-k$, $k>1$, $P(0,1)\not=0$, $P(1,0)\not=0$, and 
$P$ is monic as a polynomial in $z$ (or $\equiv 1$ if $k=d$).   
Then $g$ has a hole 
of depth $k$ at $\infty$ and no holes at 0.  Consider first the 
family, as in Example 1, given by 
 $$h_{a,t} = (atz^d + w^kP(z,w): tz^d) \in\Rat_d,$$
for $a\in\C$ and $t\in\D^*$.  Computing second iterates 
and taking a limit as $t\to 0$ gives
 $$\Phi_2(h_{a,t}) \to h_a := (aw^{kd}P(z,w)^d:w^{kd}P(z,w)^d),$$
and the degenerate $h_a$ has associated lower-degree map 
$\equiv a$.  That is, the maps $h_{a,t}$
converge to the constant $\infty$ map but their second iterates 
converge to the constant $a$.  Furthermore, $h_a\in I(d^2)$ if 
and only if $P(a,1)=0$.  By Theorem 1(a), when $P(a,1)\not=0$,
the maximal measures of $h_{a,t}$ converge weakly to 
 $$\mu_{h_a} = \frac{k}{d}\delta_\infty + 
               \frac{1}{d} \sum_{z: P(z)=0} \delta_z = \mu_g.$$
These measures do not depend on $a$.  

Let us now generalize Example 1 in the following way:  
for each $a\in\C^*$ and $t\in\D^*$, consider 
 $$g_{a,t}(z:w) = (at^kz^d + w^kP(z,w): tz^{d-k+1}w^{k-1}) \in\Rat_d.$$
As $t\to 0$, we have $g_{a,t}\to g$ in $\Ratbar_d$.  
The second iterate $\Phi_2(g_{a,t})$ has the form 
\begin{eqnarray*}
(at^k w^{kd}P^d + t^kz^{k(d-k+1)}w^{k(k-1)+k(d-k)}P^{d-k} + O(t^{k+1}): \\
       t^kz^{(d-k+1)(k-1)}w^{(k-1)^2+k(d-k+1)}P^{d-k+1} + O(t^{k+1})),
\end{eqnarray*}
so that in the limit as $t\to 0$, $\Phi_2(g_{a,t})$ converges to
  $$f_a := (w^{k(d-1)}P^{d-k}(aw^kP^k + z^{k(d-k+1)}):
            w^{k(d-1)+1}P^{d-k+1}z^{(d-k+1)(k-1)}).$$
Thus, the second iterates converge, away from the holes of $f_a$, to 
a map of degree $k(d-k+1)$ given in coordinates on $\C$ by 
  $$\phi_a(z) = \frac{z^{k(d-k+1)} + aP(z)^k}{z^{(k-1)(d-k+1)}P(z)}.$$
Since $P$ is monic of degree $d-k$, $\phi_a$ has a parabolic 
fixed point at $\infty$ for all $a\in\C^*$.  The point $\infty$ is a hole
for $f_a$ of depth $k(d-1)$, so with Lemma
\ref{pointmass} we  compute, 
 $$\mu_{f_a}(\{\infty\}) = \frac{1}{d^2} \sum_{n=0}^\infty
\frac{k(d-1)}{d^{2n}}  = \frac{k}{d+1} > 
\frac{k}{d+k},$$
and we see how the mass at $\infty$ compares with the lower bound
of Theorem 1(b).  

Finally, for different values of $a\in\C^*$, the measures $\mu_{f_a}$
are distinct.  For $k>1$, the point 0 has the same $\mu_{f_a}$-mass
for all $a\in\C^*$ because it is a preimage of the hole at $\infty$.  
However, the preimages of 0 by $\phi_a$ vary with $a$, so it is
not hard to see that the measures must vary too.  

\medskip\noindent {\bf Question.}
What is the best lower bound in the statement of Theorem 1(b) 
for $d_c >1$?  In Example 2, the limiting mass at the constant value 
$c = \infty$ is $d_c/(d+1)$.

\medskip\noindent{\bf Proof of Theorem 2}.  
The equivalence of properties (i), (ii), and (iii) was established
by Lemma \ref{iterate} in Section \ref{iteration}.  
For the implication (iv) implies (i), note first 
that within the space $\Rat_d$, the map $f\mapsto\mu_f$ is continuous
by \cite[Thm B]{Mane}.  Also, if $g\not\in I(d)$ is degenerate, 
Theorem 1(a)
implies that $f\mapsto\mu_f$ extends continuously from $\Rat_d$
to $g$.    Suppose now that 
$g_k \to g$ where $g_k\not\in I(d)$ is degenerate for all $k$
and $\mu_{g_k} \not\to \mu_g$.  Then there exists an open set
$U$ in $M^1(\P^1)$, the space of probability measures on $\P^1$,
such that $U$ contains infinitely many of the measures $\mu_{g_k}$
but $\mu_g\not\in \overline{U}$.  For each $k$ with $\mu_{g_k}\in U$,
there exists $f_k\in\Rat_d$ with $\mu_{f_k}\in U$ by Theorem 1(a).
However, in this way we can construct a sequence of non-degenerate 
rational maps 
converging to $g$ in $\Ratbar_d$ but such that $\mu_{f_k}\not\to \mu_g$.
Theorem 1(a) implies that $g\in I(d)$.  

Conversely, fix $g = H\phi\in I(d)$.  By the definition of $I(d)$, 
$\phi$ must be constant, so by a change of coordinates we can 
assume that $g(z:w) = (w^k P(z,w):0)$, where $P$ is 
homogeneous of degree $d-k$, so that $\phi$ is the
constant infinity map and $\infty$ is a hole of $g$ of depth 
$k$.  We can also assume that $P(0,1)\not=0$ so that 
0 is {\em not} a hole of $g$.  Thus, $g$ is exactly one of Examples
1 or 2 above, depending on the depth $k$ at $\infty$.  In 
each case, the family $g_{a,t}$, $a\in\C^*$, converges to $g$
as $t\to 0$ and demonstrates the discontinuity
of both the iterate map $\Phi_2$ and of $f\mapsto\mu_f$ at $g$.  
This completes the proof of Theorem 2.    
\qed

%%%%%%%%%%%%

\section{Limiting metrics}  \label{metrics}

In this section, we present details for Corollary \ref{metricprop} 
which reinterprets Theorem 1(a) in terms of conformal metrics
on the Riemann sphere:  a degenerating sequence of rational maps
has a convex polyhedral limit (with countably many vertices).  We 
apply the work of Reshetnyak on conformal metrics in planar 
domains \cite{Reshetnyak} and invoke 
the realization theorem of Alexandrov \cite[VII \S 7]{Alexandrov}
for metrics of non-negative 
curvature on a sphere to put this metric convergence into context.  

In Section \ref{dynamics} we discussed the escape rate function $G_F: 
\C^2\to\R\cup\{-\infty\}$ of a rational map $f\in\Rat_d$.  
As explained in \cite[\S 12]{D:lyap}, a hermitian metric on the 
tautological bundle $\tau\to\P^1$ is defined by
\begin{equation} \label{metric}
  \|v\|_F = \exp G_F(v)
\end{equation}
for all $v\in\C^2 -0$.  A metric $\rho_f$ is then induced on the tangent
bundle $T\P^1\iso\tau^{-2}$, uniquely up to scale.  Since $G_F$ 
is a potential function for the measure of maximal entropy $\mu_f$
of $f\in\Rat_d$, in the sense that $dd^cG_F = \pi^*\mu_f$ in $\C^2$, 
we find that the curvature form of $\rho_f$ (in the sense of
distributions)  is $4\pi\mu_f$.  In particular, the metric is Euclidean 
flat on the Fatou components.   

\medskip\noindent{\bf Example}.  Let $p\in\Rat_d$ be a polynomial
and $P(z,w) = (p(z/w)w^d,w^d)$ a lift of $p$ to $\C^2$.
The escape rate function of $p$ in $\C$ defined by 
 $$G_p(z) = \lim_{n\to\infty} \frac{1}{d^n} \log^+ |p^n(z)|$$
satisfies 
 $$G_P(z,w) = G_p(z/w) + \log|w|,$$
where $G_P$ is the escape rate function of $P$ in $\C^2$ 
\cite[Prop 8.1]{Hubbard:Papadopol}.  The associated conformal 
metric on $\C$ is given by 
  $$\rho_p = e^{-2G_p(z)}|dz|,$$
which is isometric to the flat planar metric on the filled Julia set,
  $$K(p) = \{z\in\C: G_p(z)=0\} = \{z\in\C: p^n(z) \not\to \infty 
                 \mbox{ as } n\to\infty\}.$$

\medskip
A theorem of Alexandrov says that any intrinsic metric
of non-negative curvature on $S^2$ (or $\P^1$) can be realized as the
induced metric on a convex surface in $\R^3$ (or possibly 
a doubly-covered convex planar domain).
In particular, for each $f\in \Rat_d$, the 
metric $\rho_f$ on $\P^1$ can be identified with a
convex shape in $\R^3$, unique up to scale and the isometries of 
$\R^3$.  
%Conversely, too, the induced metric on any convex surface
%in $\R^3$ can be uniformized to define a conformal metric on $\P^1$.  
See \cite[VII \S 7]{Alexandrov} or \cite[Ch.1]{Pogorelov}.
Corollary \ref{metricprop} addresses the question of how these 
metrics degenerate in an unbounded family of rational maps.  

Every probability measure $\mu$ on the Riemann sphere 
determines a metric with singularities on $\tau\to\P^1$ (and therefore on the 
tangent bundle $T\P^1$), unique up to scale, with (distributional) 
curvature equal to $4\pi\mu$.  Indeed, by 
\cite[Thm 5.9]{Fornaess:Sibony} there is a logarithmic potential
function $G_\mu$ in $\C^2$ such that $dd^c G_\mu = \pi^*\mu$,
and $G_\mu$ is unique up to an additive constant.  The metric
can be defined by equation (\ref{metric}) with $G_\mu$ in 
place of $G_F$.  That this metric induces 
a well-defined distance function on $\P^1$ follows from 
work of Reshetnyak, described in \cite[\S7]{Reshetnyak}.  

If there exists a point $z_0\in\P^1$
such that $\mu(\{z_0\}) = m >0$, then the metric is locally 
represented by 
  $$|z-z_0|^{-2m} e^{u(z)} |dz|,$$
where $u$ is subharmonic near $z_0$.   For $m\geq 1/2$, 
the point $z_0$ will be at infinite distance from all other points
in the sphere, and $z_0$ will be called an {\bf infinite end} of this metric.  
There can be at most two such points.  

%A detailed description of the relation between conformal metrics
%in $\C$ and logarithmic potentials and curvature distributions is 
%given in \cite[\S 7]{Reshetnyak}.  For example, suppose a metric has 
%curvature distribution $4\pi\mu$ for a probability measure $\mu$ with 
%an isolated atom at $z_0\in\C$, so that $m=4\pi\mu(z_0) >0$.
%Then the metric is locally represented by 
%This formulation shows that a neighborhood of $z_0$ is isometric
%to a Euclidean cone of angle $2\pi - m$.   If the cone angle is 
%non-positive, then the point $z_0$ is at infinite distance from all 
%points, so that $z_0$ is an {\bf infinite end} of this metric.

We will say that a metric on $\P^1$ is {\bf convex polyhedral} if its 
curvature distribution can be expressed as a countable sum 
of delta masses.  For example, the curvature measure
of the induced metric on a convex
polyhedron with finitely many vertices in $\R^3$
is a finite sum of delta masses.  
The degenerate
maps in $\Ratbar_d$ have associated metrics which 
are convex polyhedral, by the definition of the probability measures $\mu_f$.

Reshetnyak proved the following convergence theorem about 
these conformal metrics \cite[Thm 7.3.1]{Reshetnyak}.  Suppose 
that $\rho_k$ and $\rho$ are metrics on $\P^1$ with curvature
distributions $4\pi\mu_k$ and $4\pi\mu$, where $\mu_k$ and 
$\mu$ are probability measures.  
If $\mu_k\to\mu$ weakly, then the metrics $\rho_k$ (as distance
functions on $\P^1\times\P^1$) converge to $\rho$, locally uniformly 
on the complement of any points $z\in\P^1$ with $\mu\{z\}\geq 1/2$.  
That is to say, the convergence is uniform away from any infinite ends 
of the metric $\rho$.  

Let it be noted that such a convergence theorem 
requires a normalization, a choice of scale for these metrics.  If, for 
example, points $0$ and $\infty$ in $\P^1$ are not infinite ends for any
of the metrics, we could fix $\rho_k(0,\infty) = \rho(0,\infty) = 1$
for all $k$.  

\medskip\noindent {\bf Proof of Corollary \ref{metricprop}}.  
Suppose the sequence $f_k\in\Rat_d$ converges in $\Ratbar_d$
to $f\not\in I(d)$.  The metric associated to $f$ is convex polyhedral since
$\mu_f$ is a countable sum of delta masses.  Choose points $a, b
\in\P^1$ which are not infinite ends for the metric of $f$.  Let 
$\rho_k$ and $\rho$ denote the distance functions on $\P^1\times\P^1$
defined by $\mu_{f_k}$ and $\mu_f$ respectively, normalized so
that $\rho_k(a,b) = \rho(a,b) = 1$ for all $k$.  By Theorem 1(a), 
the measures $\mu_{f_k}$ converge weakly to $\mu_f$.  Therefore,
by \cite[Thm 7.3.1]{Reshetnyak}, the metrics $\rho_k$ converge
to $\rho$ locally uniformly on the complement of the infinite
ends of $\rho$.  
\qed

\medskip\noindent
{\bf Question}.  Are {\em all} limiting metrics in the boundary of $\Rat_d$ 
polyhedral?  That is, if $f_k\to f\in I(d)$ such that the maximal measures
$\mu_{f_k}$ converge weakly, is the limiting distribution a countable 
sum of delta masses?   

\medskip\noindent
{\bf Polynomial limits.}  In the next section, 
we will see from Proposition \ref{polylimit}
that all measures of the form 
  $$\mu = \frac{1}{d} \sum_i \delta_{z_i},$$
where $\{z_i: i=1,\ldots,d\}$ is any collection of $d$ (not necessarily 
distinct) points in $\P^1$,
arise as limits of the maximal measures for degree $d$
polynomials.    Metrically, these 
limits correspond to all convex 
polyhedra with $d$ vertices of equal cone angle, 
the objects of study in \cite{Thurston:shapes}.  In our case, several 
or all vertices may coalesce.  
When the limit measure is $\delta_\infty$, for example, the metric is that 
of the flat plane.  

%\medskip\noindent
%{\bf Question.}   Are all limiting metrics in the boundary of the space
%of polynomials,
%$\Poly_d$, polyhedral with {\em finitely many} vertices?  

%%%%%%%%%%%

\section{Further examples}  \label{examples}

In this final section, we study the above ideas as they apply
to the boundary of the space of polynomials 
$\Poly_d \subset \Rat_d$ in $\Ratbar_d\iso\P^{2d+1}$.  
We also explain how Epstein's sequences of degree 2 rational
maps in \cite{Epstein:bounded} arise and achieve the lower bound
of Theorem 1(b).  

\medskip\noindent
{\bf The boundary of $\Poly_d$}.  
The space of polynomials of degree $d\geq 2$ satisfies $\Poly_d
\iso \C^*\times\C^d$ and $\Polybar_d \iso \P^{d+1}$ in $\Ratbar_d$.
The boundary $\del\Poly_d$ has two irreducible components.  
Indeed, a point $p\in\Poly_d$ can be expressed in homogeneous
coordinates by 
  $$p(z:w) = (a_dz^d + a_{d-1}z^{d-1}w + \cdots + a_0w^d: b_0w^d),$$
so that $\del\Poly_d$ is defined by $\{a_db_0 = 0\}$.  Furthermore, 
  $$\Polybar_d \cap I(d) = \{a_d = b_0 = 0 \}$$
is the intersection of the two boundary components, and it consists
of the ``constant $\infty$" polynomials with a hole at $\infty$.  It has 
codimension 2 in $\Polybar_d$.  

If a sequence of polynomials converges locally
uniformly in $\C$ to a polynomial of lower degree, then 
we have $a_d=0$ and $b_0\not= 0$ in the limit.  For all 
points $p$ in
the locus $\{a_d = 0 \mbox{ and } b_0\not=0\}$, the associated
measure is $\mu_p = \delta_\infty$:  the only hole is at $\infty$.
By Theorem 1(a), the measures of maximal entropy will 
converge to $\delta_\infty$.  However, 
the supports of the measures, the Julia sets of the polynomials, 
do not necessarily
go off to infinity (in the Hausdorff topology).  

\medskip\noindent{\bf Example.}
Consider the family of cubic polynomials 
  $$p_\eps(z) =  \eps z^3 + z^2$$ 
for $\eps\in\C^*$.  As $\eps \to 0$, we have $p_\eps \to (z^2w: w^2)$
in $\Polybar_3$.  As maps, $p_\eps \to z^2$ locally uniformly in 
$\C$ by Lemma \ref{convergence}, but by Theorem 1(a), 
$\mu_{p_\eps} \to \delta_\infty$ weakly.  
For small $\eps$, $p_\eps$ is polynomial-like 
of degree 2 in a neighborhood of the unit disk.  
The Julia set of the 
polynomial-like restriction of $p_\eps$ is part of the Julia set for $p_\eps$.
It converges to the unit
circle as $\eps\to 0$, but carries no measure.  (From the point of view
of external rays, almost none land so deeply.)  The components of the 
Julia set which carry the measure 
are tending to infinity as $\eps \to 0$.  

\medskip
In general, we have the following corollary of Theorem 1(a).  

\begin{cor}  \label{genericpoly}
Let $p_k$ be a sequence of polynomials in $\Poly_d$ such that 
$p_k \to p \not\in I(d)$ in $\Polybar_d$.  Then the 
measures $\mu_{p_k}$ converge weakly to a measure 
$\mu_p$ supported in at most $d$ points.  
\end{cor}

\proof
For notational simplicity, let $A = \{a_d =0\}$ and $B = \{b_0 =0\}$ be 
the irreducible components of $\del\Poly_d$.  For each $p = H\phi_p 
\in A- I(d)$, we have $\mu_p = \delta_\infty$, since $\phi_p$ is a 
polynomial of degree $< d$ and the only hole is at $\infty$.  For 
$p = H\phi_p \in B - I(d)$, $\phi_p$ is the constant $\infty$, the
$d$ holes are finite, and $\mu_p$ is supported at the finite holes.  
By Theorem 1(a),  $\mu_{p_k} \to \mu_p$ 
weakly as $k\to\infty$.  
\qed

\begin{prop} \label{polylimit}  For $d\geq 2$,
let $z_1, z_2, \ldots, z_d$ be (not necessarily distinct) points
in $\P^1$, and let 
  $$ \mu = \frac{1}{d} \sum_{i=1}^d \delta_{z_i}$$
be the probability measure supported equally at these points.  There
exists a sequence of polynomials $p_k \in \Poly_d$ such that 
$\mu_{p_k} \to \mu$ weakly.  
\end{prop}

\proof
Suppose first that no $z_i$ is $\infty$.  Let $P(z,w)$ be a homogeneous
polynomial of degree $d$ with roots in $\P^1$ at the points $z_i$.  
Consider the sequence
$$p_k(z:w) = (P(z,w): w^d/k)\in\Poly_d$$
as $k\to\infty$.  The limit $f = (P(z,w):0)$ is the constant 
$\infty$ map with holes at the roots of $P$, so in particular $f\not\in I(d)$.
By Theorem 1(a), the measures of maximal entropy for $p_k$ converge weakly to 
$\mu_f = \sum_{P(z,1)=0} \delta_z/d$.  

For an arbitrary collection of $d$ points $\{z_i\}$, the given measure
$\mu$ can be approximated by measures of the form of $\mu_f$ described
above.  Therefore, there must exist a sequence of polynomials with 
maximal measures limiting on $\mu$.  
\qed

\begin{cor} \label{denseimage}
Given any probability measure $\mu$ on $\P^1$, there is a 
sequence of polynomials $p_k$, of degrees tending to infinity,
such that $\mu_{p_k} \to \mu$ weakly.
\end{cor}

\proof
The finite atomic measures with rational mass at every point 
are dense in the space of probability measures.  By Proposition 
\ref{polylimit}, each such 
measure is a limit of measures of maximal entropy in 
$\Poly_d$ for a sufficiently
large degree $d$.  
\qed

\medskip
%{\bf The indeterminacy locus at the boundary of $\Poly_d$}.
%It is straightforward to construct a sequence of polynomials in $\Poly_d$
%converging to $I(d)$ for which some iterate has a positive degree
%limit, simply by forcing a critical point to be periodic.  
%For example, consider the family 
%  $$q_\eps(z) = \eps^3 z^3 + \eps z^2 + c(\eps),$$
%where the constant $c(\eps)$ is chosen so that the critical 
%point at 0 has period two and $c(\eps) = O(1/\eps^2)$ for $\eps$ 
%small.  (This $c(\eps)$ is obtained by choosing the ``correct'' 
%square root.)  Since the constant term is unbounded, $q_\eps$ 
%converges to the constant $\infty$ as $\eps\to 0$, locally uniformly
%on $\C$.  Algebraically, $q_\eps \to (w^3:0)\in I(3)$ in $\Polybar_3$,
%so that $\infty$ is a hole of depth 3.
%Computing the second iterate $\Phi_2(q_\eps)$, we see that as $\eps\to 0$, 
%  $$\Phi_2(q_\eps) \to z^2$$
%locally uniformly on $\C$.   In this example, because the second
%iterate does not limit in $I(9)$, we can compute the limiting measures
%and find that $\mu_{q_\eps} \to \delta_\infty$.  

It should be noted at this point that, while the limiting measures 
away from $I(d)$ are supported in at most $d$ points, there is no
bound on the number of points in the support of a general 
limit at the boundary of $\Poly_d$.
For example, various normalizations of the family $\eps z^3 + z^2$
as $\eps\to 0$ can give a limiting measure with 1, 2,  or $2+2^n$ points 
in its support for any desired $n\geq 1$.  

\medskip\noindent {\bf Question.}
What is the closure of $\Poly_d$ in the space of probability measures?

%\medskip\noindent
%{\bf General limits at boundary of $\Poly_d$.}
%From experimental evidence, it appears that 
%even when a sequence of polynomials converges to $I(d)$, the limiting 
%measure should be supported at finitely many points.    In degree 
%two, it is possible to show (using the fact that 
%the moduli of certain 
%annuli separating Green level lines grow unboundedly if the sequence
%diverges in moduli space) that the limiting measures can 
%be supported in at most three points.
% and have the form 
%  $$\frac{1}{2^n} \delta_a + \frac{1}{2^n} \delta_b 
%    +\left( 1 - \frac{1}{2^{n-1}} \right)  \delta_\infty,$$
%for some $n\geq 1$ and any pair of points $a$ and $b$ in $\P^1$.   
%Different limiting measures arise
%from different normalizations.  For example, choosing a fixed 
%point and a particular $n$-th preimage of that fixed point to remain
%finite, we obtain points in the plane of mass $1/2^n$ with the 
%remaining mass at $\infty$.  

\medskip\noindent
{\bf Epstein's sequences in $\Rat_2$.}  In \cite[Prop 2]{Epstein:bounded}, 
Epstein gave 
the first example of the discontinuity of the iterate map at 
the boundary of the space of rational maps.  He 
studied unbounded sequences in the moduli 
space of degree 2 rational maps, and in particular, examined sequences 
of rational maps in $\Rat_2$ converging to a 
degenerate map $f=H\phi\in \Ratbar_2$ for 
which $\phi$ is an
elliptic Mobius transformation of order $q>1$.  If a sequence
of rational maps 
approaches $f$ from a particular direction in
$\Ratbar_2 \iso \P^5$ 
(depending on a complex
parameter $T$), there are certain conjugates of this sequence 
converging to $I(2)$ such that 
their $q$-th iterates converge to the
degree 2 map, $\phi_T(z) = z + T + 1/z$ (uniformly away from the
holes at 0 and $\infty$).  

When $q=2$, his examples realize the lower bound in Theorem 1(b).   
Indeed, he provides a family of sequences $F_{T,k}\in\Rat_2$, 
$T\in\C$,
normalized
so that the critical points of $F_{T,k}$ are at 1 and $-1$ for all 
$T\in\C$ and all $k\geq 1$, and 
such that
$F_{T,k} \to (zw:0) \in I(2)$ as $k\to \infty$.  
For all $T\in\C$, this limit is 
the constant infinity map with holes at 0 
and $\infty$, each of depth 1.  The second iterates of $F_{T,k}$ 
converge to $F_T = (zw(z^2 + Tzw + w^2):z^2w^2) \not\in I(4)$ 
as $k\to\infty$.  By Theorem 1(a), 
the measures $\mu_{F_{T,k}}$ converge weakly to $\mu_{F_T}$ 
as $k\to\infty$, and 
we can compute with Lemma \ref{pointmass} that 
  $$\mu_{F_T}(\{\infty\}) = \frac{1}{4} \sum_{n=0}^\infty \frac{1}{4^n}
                                     = \frac{1}{3},$$
since $F_T$ has a hole of depth 1 at $\infty$ which is fixed by $\phi_T$.
Recalling that the limit of the first iterates was in $I(2)$ and had a hole
of depth 1 at $\infty$, we see that 
the limiting measure at $\infty$ is exactly the lower bound in 
Theorem 1(b) when the degree is 2 and the depth is 1.  

Epstein used these examples in his proof that certain hyperbolic
components in the moduli space $\Rat_2/\PSL_2\C$ are bounded.  
It is my hope that a more systematic understanding of the 
iterate map and the boundary of the space of rational maps
will be applicable to related questions in general degrees.

§
\bigskip
\def\cprime{$'$}

%\bibliographystyle{../tex/bib/math}
%\bibliography{../tex/bib/math}

\begin{thebibliography}{FLM}

\bibitem[Al]{Alexandrov}
A.~D. Alexandrov.
\newblock {\em Vnutrennyaya {G}eometriya {V}ypuklyh {P}overhnoste\u\i}.
\newblock OGIZ, Moscow-Leningrad, 1948.

\bibitem[De1]{D:lyap}
L.~DeMarco.
\newblock {Dynamics of rational maps: {L}yapunov exponents, bifurcations, and
  capacity}.
\newblock {\em Math. Ann.} {\bf 326}(2003), 43--73.

\bibitem[De2]{D:moduli2}
L.~DeMarco.
\newblock {The boundary of the moduli space of quadratic rational maps}.
\newblock {\em Preprint, \em{2004}}.

\bibitem[Ep]{Epstein:bounded}
A.~Epstein.
\newblock {Bounded hyperbolic components of quadratic rational maps}.
\newblock {\em Ergodic Theory Dynam. Systems} {\bf 20}(2000), 727--748.

\bibitem[FS]{Fornaess:Sibony}
J.~E. Forn{\ae}ss and N.~Sibony.
\newblock {Complex dynamics in higher dimensions}.
\newblock In {\em Complex Potential Theory (Montreal, PQ, 1993)}, pages
  131--186. Kluwer Acad. Publ., Dordrecht, 1994.

\bibitem[FLM]{FLM}
A.~Freire, A.~Lopes, and R.~Ma\~{n}{\'e}.
\newblock {An invariant measure for rational maps}.
\newblock {\em Bol. Soc. Brasil. Mat.} {\bf 14}(1983), 45--62.

\bibitem[HP]{Hubbard:Papadopol}
J.~Hubbard and P.~Papadopol.
\newblock {Superattractive fixed points in ${\bf C}^n$}.
\newblock {\em Indiana Univ. Math. J.} {\bf 43}(1994), 321--365.

\bibitem[Ly]{Lyubich:entropy}
M.~Lyubich.
\newblock {Entropy properties of rational endomorphisms of the {R}iemann
  sphere}.
\newblock {\em Ergodic Theory Dynamical Systems} {\bf 3}(1983), 351--385.

\bibitem[Ma1]{Mane:unique}
R.~Ma\~{n}{\'e}.
\newblock {On the uniqueness of the maximizing measure for rational maps}.
\newblock {\em Bol. Soc. Brasil. Mat.} {\bf 14}(1983), 27--43.

\bibitem[Ma2]{Mane}
R.~Ma\~{n}{\'e}.
\newblock {The {H}ausdorff dimension of invariant probabilities of rational
  maps}.
\newblock In {\em Dynamical Systems, Valparaiso 1986}, pages 86--117. Springer,
  Berlin, 1988.

\bibitem[Mi]{Milnor:quad}
J.~Milnor.
\newblock {Geometry and dynamics of quadratic rational maps}.
\newblock {\em Experiment. Math.} {\bf 2}(1993), 37--83.
\newblock With an appendix by the author and Lei Tan.

\bibitem[Po]{Pogorelov}
A.~V. Pogorelov.
\newblock {\em Extrinsic geometry of convex surfaces}.
\newblock American Mathematical Society, Providence, R.I., 1973.
\newblock Translated from the Russian by Israel Program for Scientific
  Translations, Translations of Mathematical Monographs, Vol. 35.

\bibitem[Re]{Reshetnyak}
Yu.~G. Reshetnyak.
\newblock {\em Geometry. {IV}}, volume~70 of {\em Encyclopaedia of Mathematical
  Sciences}.
\newblock Springer-Verlag, Berlin, 1993.
\newblock Nonregular Riemannian geometry, A translation of {\it Geometry, 4
  (Russian)}, Akad.\ Nauk SSSR, Vsesoyuz.\ Inst.\ Nauchn.\ i Tekhn.\ Inform.,
  Moscow, 1989 [MR 91k:53003], Translation by E. Primrose.

\bibitem[Si]{Silverman}
J.~H. Silverman.
\newblock {The space of rational maps on $\bf {P}\sp 1$}.
\newblock {\em Duke Math. J.} {\bf 94}(1998), 41--77.

\bibitem[Th]{Thurston:shapes}
W.~P. Thurston.
\newblock {Shapes of polyhedra and triangulations of the sphere}.
\newblock In {\em The Epstein birthday schrift}, volume~1 of {\em Geom. Topol.
  Monogr.}, pages 511--549 (electronic). Geom. Topol. Publ., Coventry, 1998.

\end{thebibliography}

\end{document}